\documentstyle[12pt] {article}

\def\zid{1\kern-0.36em\llap~1}

\newcommand{\beq}{\begin{equation}}
\newcommand{\ber}{\begin{eqnarray}}
\newcommand{\eeq}{\end{equation}}
\newcommand{\eer}{\end{eqnarray}}

\begin{document}

\rightline{ SUNY BING 6/25/96r }
\vspace{1mm}
\begin{center}
{\large \bf ON THE TWO q-ANALOGUE \\
LOGARITHMIC FUNCTIONS: $\ln_q(w)$, $\ln \{ e_q(z) \}$
}\\[2mm]
Charles A. Nelson\footnote{Electronic address: cnelson @
bingvmb.cc.binghamton.edu . } and Michael G.
Gartley\footnote{Present address:  Rochester Institute of
Technology, 1 Lomb Memorial Drive, Rochester, N.Y. 14623.}
\\
{\it Department of Physics, State University of New York at
Binghamton\\
Binghamton, N.Y. 13902-6016}\\[2mm]
\end{center}


\begin{abstract}
There is a simple, multi-sheet Riemann surface associated
with $e_q(z)$'s inverse function $ln_q(w)$ for $0<q \leq 1$.
A
principal sheet
for $ln_q(w)$ can be defined. However, the topology of the
Riemann surface for $ln_q(w)$ changes each time $q$ increases
above the collision point $q_{\tau}^*$ of a pair of the
turning points $\tau_i$ of $e_q(x)$.  There is also a power
series
representation for $ln_q(1+w)$.  An infinite-product
representation for $e_q(z)$ is used to obtain the ordinary
natural logarithm $\ln \{e_q(z) \}$ and the values of the sum
rules $\sigma _n^e\equiv
\sum\limits_{i=1}^\infty
\left( \frac 1{z_i}\right) ^n $ for the zeros $z_i$ of
$e_q(z)$.
For $|z|<|z_1|$, $e_q(z)=exp\{ b(z) \}$ where $ b(z)=-
\sum\limits_{n=1}^\infty \frac 1n\sigma _n^e z^n $.  The
values of the sum rules for the q-trigonometric functions,
$\sigma _{2n}^c$ and $\sigma _{2n+1}^s$, are q-deformations
of the usual Bernoulli numbers.
\end{abstract}

\newpage

\section{Introduction:}

The ordinary exponential and logarithmic functions find
frequent and varied applications in all fields of physics.
Recently in the study of
quantum algebras, the $q$-exponential function [1] or mapping
$w=e_q(z)$ has reappeared [2-4]
from a rather dormant status in mathematical physics.  This
order-zero entire function can be defined by
\begin{equation}
e_q(z)\equiv \sum_{n=0}^\infty \frac{z^n}{[n]!}
\end{equation}
where
\begin{equation}
[n]=\frac{q^{n/2}-q^{-n/2}}{q^{1/2}-q^{-1/2}}
\end{equation}
The series in Eq.(1) converges uniformly and absolutely for
all finite $z$. Since $[n]$ is invariant under $q \rightarrow
1/q $, for real $q$
it suffices to study $0<q \leq 1$. The q-factorial is defined
by $[n]!\equiv
[n][n-1]\cdots [1]$, $[0]!\equiv 1$.  As $q \rightarrow 1$,
$e_q(z) \rightarrow exp(z)$ the ordinary exponential
function.

In [5], we reported some of the remarkable analytic and
numerical properties of the infinity of zeros, $z_i$, of
$e_q(x)$ for $x<0$. In particular, as $q$ increases above the
first
collision point at $q_z^* \approx 0.14 $, these zeros collide
in
pairs and then move off into the complex $z$ plane, see Fig.
1.
They move off as (and remain) a complex conjugate pair
$\mu_{A,\bar{A} }$. The
turning points of $e_q(z)$, i.e. the zeros of the first
derivative $e_q^{^{\prime }}(z)\equiv de_q(z)/dx$, behave in
a similar
manner. For instance, at $q_{\tau}^* \approx 0.25 $ the first
two turning points,
$\tau_1$ and $\tau_2$, collide and move off as a complex
conjugate pair $\tau_{A,\bar{A} }$.

In this paper, we first show that there
is a simple, multi-sheet
Riemann surface associated with $w=e_q(z)$'s inverse function
$z=ln_q(w)$.  As with the usual $ln(w)$ function, the
Riemann surface of $z=ln_q(w)$ defines a single-valued map
onto the entire complex $z$ plane. Also, as in the usual case
when $q=1$, a principal sheet for
$z=ln_q(w)$ can be defined.  However, unlike for the ordinary
$ln(w)$ and
$exp(z)$, the topology of the Riemann surface for $ln_q(w)$
changes each time $q$ increases above the collision point
$q_{\tau}^*$
of a pair of the turning points $\tau_i$ of $e_q(z)$. The
turning points of $e_q(z)$ can be used to
define square-root
branch points of $ln_q(w)$ in the complex $w$ plane, i.e.
$b_i = e_q(\tau_i)$.

In Sec. 3, we obtain a power series representation for
$ln_q(1+w)$.

In the mathematics and physics literature\footnote{Recent
reviews of quantum algebras are listed in [6].}, one
also finds the
exponential function $E_q(z)$ defined by Jackson[7-8].  It
also
is
given by Eq.(1) but with $[n]$ replaced by $[n]_J$
where
\begin{equation}
[n]_J=q^{(n-1)/2}[n]=\frac{1-q^n}{1-q}
\end{equation}
For $q>1$, $E_q(z)$ has simpler properties\footnote{ For
$0<q<1$, $E_q(z)$ is a meromorphic function whose power
series converges
uniformly
and absolutely for $|z|<(1-q)^{-1}$ but diverges otherwise.
However by the
relation, $E_s(x)E_{1/s}(-x)=1$ for $s$ real, results for
$q>1$ can be used
for $0<q<1$, see Ref. [5]. }than $e_q(z)$.  We
also construct the Riemann
surface for its inverse function $Ln_q(w)$.  With the
substitution $[n] \rightarrow
[n]_J $, the power series
representation for $ln_q(1+w)$ also holds for $Ln_q(1+w)$.

Second, in Sec. 4, we use the infinite-product representation
[5] for $e_q(z)$ to (i) obtain the ordinary natural logarithm
$\ln \{ e_q(z) \}$, and to (ii) evaluate for arbitrary
integer
$n
>
0$ the sum rules
\begin{equation}
\sigma _n^e\equiv \sum_{i=1}^\infty \left( \frac
1{z_i}\right)
^n
\end{equation}
for the zeros $z_i$ of $e_q(z)$. Therefore, for c-number
arguments
\begin{equation}
e_q(x)e_q(y)=\exp \left\{ b(x)+b(y)\right\}
\end{equation}
where $b(x)$ is defined below in Eq.(20). For $|z|<|z_1|$ the
modulus of the
first zero,
\begin{equation}
b(z)=-\sum_{n=1}^\infty \frac 1n\sigma _n^ez^n
\end{equation}

We also obtain the logarithms and values of the associated
sum rules for all derivatives and integrals of
$e_q(x)$, and for the associated q-trigonometric functions
[1,5]
$cos_q (z)$ and $sin_q (z)$.  These results also hold for the
analogous functions involving $[n]_J$.

Sec. 5 contains some concluding remarks.  In particular, the
values of the sum rules for the q-trigonometric functions,
$\sigma _{2n}^c$ and $\sigma _{2n+1}^s$, are q-deformations
of the usual Bernoulli numbers.

\section{Riemann Surfaces of q-Analogue Logarithmic \newline
Functions
$ln_q(w)$ and $Ln_q(w)$:}

For two reasons, we begin by first analyzing the Riemann
surface
associated with the mapping of Jackson's exponential function
$w= u+i v = E_q(z)$ and of its inverse $z= x+i y = Ln_q(w)$.
First, the generic structure of the Riemann surface for
$Ln_q(w)$ for $q^E >1$ is the same as that for $ln_q(w)$ for
$ q^e < (q^* \approx 0.14)$.  Second, as $q^e$ varies the
topology
of
the Riemann surface changes for $ln_q(w)$ but the topology
remains invariant for $Ln_q(w)$ for all $q^E>1$.  Normally we
will suppress the superscripts ``$E$ or $e$" on the $q$'s for
there should be no confusion.

\subsection{Riemann surface for $Ln_q(w)$:}

Figs. 2 and 3 show the Riemann sheet structure and the
mappings
of Jackson's exponential function $w=E_q(z)$ and of its
inverse
$z= Ln_q(w)$ for $q^E \approx 1.09$. These figures suffice
for
illustrating the Riemann sheet for all $q>1$ because the
zeros and
turning points of $E_q(z)$ do not collide, but simply move
along
the negative $x$ axis and out to infinity as $q \rightarrow
1$.

These figures also illustrate the Riemann surface for
$w=e_q(z)$
and $z=ln_q(w)$ but only prior to the collision of the first
pair of zeros at $q \approx 0.14$.

Notice that the imaginary part $Im\{e_q(z) \} =0$ on all
``solid"
contour lines in Fig. 2b whereas the real part $Re\{e_q(z) \}
=0$
on all ``dashed" contour lines. The turning points in the
complex
$z$ plane are denoted by small dark squares, whereas their
associated branch points in $w$ are denoted by small dark
circles.

Numerically, for $q^E \approx 1.09$, the first 4 zeros of
$E_q(z)$ are located at $-12.1111, -13.2011, \newline
-14.3892, -15.6842$.  The first 4 turning points and $
Ln_q(w)$'s branch
points ($b_i$ in $10^{-11}$ units) are respectively at
$(\tau_i,b_i)= (-12.4,-43), (-13.6,5.0), (-14.9,-1.8), (-
16.3,4.4)$.  Since
$q^E \approx 1$, the asymptotic formula in
[5] for $\tau_i^E$ is a bad approximation for these values.

Figures for the lower-sheets of a Riemann surface $w$ are
omitted in this paper since they simply have the conjugate
structures, per the Schwarz reflection
principle.

\subsection{Riemann surface for $ln_q(w)$:}

For $q< \approx 0.14$, Figs. 1-3 also show the topology and
branch point structure for the mappings $w=e_q(z)$ and its
inverse $z= ln_q(w)$.

Figs. 4-5 are for after the collision of the first pair of
zeros of $e_q(z)$ but prior to the collision of the first
pair of its turning points, so the structure shown is generic
for $0.14< q < 0.25$. Note that $w_A=e_q (\mu_A) = 0$
occurs as an analytic point for $w=e_q(z)$ which is not
possible for the ordinary $exp(z)$ in the finite $z$ plane.

Numerically, Figs. 4-5 are for $q \approx 0.22$; the first 2
zeros of $e_q(z)$ are located at $\mu_A=-2.51 + i 0.87,
\mu_{\bar{A}} =\bar{\mu_A}$. The first 2 turning points and $
ln_q(w)$'s branch points ($b_i$ in $10^{-3}$ units) are
respectively at $(\tau_i,b_i)= (-2.6,47.70), (-4.7,69.36)$.

Figs. 6-8 are for after the collision of the first pair of
turning points of $e_q(z)$.  The topology of the Riemann
surface has a new inter-surface structure due to this
collision; the figures and their captions explain this new
structure. In particular versus Fig. 5, following
the collision at $q_{\tau}^* \approx 0.25$, there no longer
exists the $b_1-b_2$ passage from the lower-half of the
principal $w$ sheet to the first lower $w$ sheet.  Instead,
the $b_A$ passages are to the second upper $w$ sheet.

Numerically, Figs. 6-8 are for $q \approx 0.35$.  The first 2
zeros of $e_q(z)$ are now located at $\mu_A=-2.8222 + i
1.969, \mu_{\bar{A}} =\bar{\mu_A}$; the third zero remains on
the negative real axis at $\mu_3=-5.19755$. The first 4
turning points and $ ln_q(w)$'s branch points ($b_i$ in
$10^{-3}$
units) are respectively at $(\tau_i, b_i)= (-3.5434 \pm i
1.32945,22.2415 \pm i 18.79), (-6.3471,-9.09587),
(-10.7028,87.536)$.  In Figs. 7-8, for clarity of
illustration, the position of $b_A$ has been displaced from
its true position.

\section{Power Series Representations for $\ln_q(1+w)$
\newline
and
$Ln_q(1+w)$:}

To obtain the power series for $\ln _q(1+w)$, we write
\begin{equation}
\begin{array}{c}
\ln _q(1+w)=c_1w+c_2w^2+... \\
=\sum\limits_{n=1}^\infty c_nw^n
\end{array}
\end{equation}

Then for $a=\ln _q(1+w)$,
\begin{equation}
\begin{array}{c}
e_q^a=1+a+
\frac{a^2}{[2]!}+... \\ =1+w
\end{array}
\end{equation}

So by equating coefficients, we find
\begin{equation}
\begin{array}{c}
c_1=1 \\
c_n=-\sum\limits_{l=2}^n\frac 1{[l]!}\left\{
\sum\limits_{\left(
k_1,k_2,\cdots k_l\right)
}c_{k_1}c_{k_2}\cdots c_{k_l}\right\} ,n\geq 2
\end{array}
\end{equation}

In order to follow later expressions in this paper, it is
essential to
understand the second summation $\sum\limits_{\left(
k_1,k_2,\cdots k_l\right)}$:

In it, each $k_i=$ ``positive integer'', $i=1,2,\ldots l.$

$\left( k_1,k_2,\cdots k_l\right) $ denotes that, for fixed
$n$ and $l$,  the
summation is the symmetric permutations of each
partition of n
which satisfy the condition $``k_1+k_2+\cdots k_l=n".$

For instance, for $n=4$:
\begin{equation}
\begin{array}{c}
\sum\limits_{\left( k_1,k_2,k_3,k_4\right)
}c_{k_1}c_{k_2}c_{k_3}c_{k_4}=\{c_1c_1c_1c_1\}=(c_1)^4 \\
\sum\limits_{\left( k_1,k_2,k_3\right)
}c_{k_1}c_{k_2}c_{k_3}=\{c_1c_1c_2+c_1c_2c_1+c_2c_1c_1\}=3c_1
c_1
c_2 \\
\sum\limits_{\left( k_1,k_2\right)
}c_{k_1}c_{k_2}=\{c_2c_2\}+\{c_1c_3+c_3c_1\}=(c_2)^2+2c_1c_3
\end{array}
\end{equation}

This power series for $\ln _q(1+w)$ is expected to converge
only
for some $w$ domain, e.g. for
$w\leq $ ``modulus of distance to the nearest branch point''.
Note that
as $q \rightarrow 0$, $w=e_q(z) \rightarrow w=1+z$ and
$z=\ln_q(w) \rightarrow z=w-1$, so $e_q \{ \ln_q(w) \}
\rightarrow e_q\{ w-1 \} \rightarrow w$.

Thus, the first few terms give
\begin{equation}
\begin{array}{c}
\begin{array}{c}
\ln _q(1+w)=w-\frac 1{[2]!}w^2-\left\{ \frac 1{[3]!}-\frac
2{[2]![2]!}\right\} w^3 \\
-\left\{ \frac 1{[4]!}-\frac 2{[2]!}\left( \frac 1{[3]!}-
\frac
2{[2]![2]!}\right) +\left( \frac 1{[2]!}\right) ^3-\frac
3{[3]![2]!}\right\}
w^4+\ldots
\end{array}
\\
=w-\frac 1{[2]!}w^2-\left\{ \frac 1{[3]!}-2\left( \frac
1{[2]!}\right)
^2\right\} w^3 \\
-\left\{ \frac 1{[4]!}-\frac 5{[3]![2]!}+5\left( \frac
1{[2]!}\right)
^3\right\} w^4+\ldots
\end{array}
\end{equation}

Notice that here the q-derivative operation defines a new
function, $d\ln _q(w) / {d_qw} \equiv \ln_q(w)^{'}  \neq
\frac
1{w}$,
because it does {\em not} yield a known q-special function
since
\begin{equation}
\begin{array}{c}
\frac d{d_qw}\ln _q(1+w)=1-w-\left\{ \frac 1{[2]!}-2[3]\left(
\frac
1{[2]!}\right) ^2\right\} w^2 \\
-\left\{ \frac 1{[3]!}-\frac{5[4]}{[3]![2]!}+5[4]\left( \frac
1{[2]!}\right)
^3\right\} w^3+\cdots
\end{array}
\end{equation}
unlike [5] for $e_q(z)$, $\cos _q(z)$, and $\sin
_q(z)$.

\section{Natural Logarithms and Sum Rules for $e_q(z)$
\newline and Related Functions:}

By the Hadamard-Weierstrass theorem, it was shown in Ref.[5]
that the following order-zero entire functions have infinite
product representations in terms of their respective zeros:
\begin{equation}
e_q(z)=\prod_{i=1}^\infty \left( 1-\frac z{z_i}\right)
\end{equation}
\begin{equation}
\begin{array}{c}
e_q^{(r)}(x)\equiv
\frac{d^r}{dx^r}e_q(x)=\alpha _r\prod\limits_{i=1}^\infty
\left( 1-
\frac
x{z_i^{(r)}}\right) ;r=1,2,\ldots  \\ \alpha
_r=\frac{r!}{[r]!}
\end{array}
\end{equation}
\begin{equation}
\begin{array}{c}
e_q^{(-r)}(x)=\int^xdx_1\int^{x_1}dx_2\ldots
\int^{x_r}dx_re_q(x_r)+
{poly. deg. (r-1)},r\geq 1 \\
\equiv
\sum\limits_{n=0}^\infty
\frac{n!}{(n+r)!}\frac{x^{n+r}}{[n]!} \\ =\left(
\frac{x^r}{r!}\right)
\prod\limits_{i=1}^\infty \left( 1-\frac x{z_i^{(-r)}}\right)
\end{array}
\end{equation}
\begin{equation}
\begin{array}{c}
\cos _q(z)\equiv \sum\limits_{n=0}^\infty (-)^n
\frac{z^{2n}}{[2n]!} \\ =\prod\limits_{i=1}^\infty \left( 1-
\left(
\frac
z{c_i}\right) ^2\right)
\end{array}
\end{equation}
\begin{equation}
\begin{array}{c}
\sin _q(z)\equiv \sum\limits_{n=0}^\infty (-)^n
\frac{z^{2n+1}}{[2n+1]!} \\ =z\prod\limits_{i=1}^\infty
\left( 1-
\left(
\frac
z{s_i}\right) ^2\right)
\end{array}
\end{equation}

\subsection{Derivation of $\ln \{e_q(z) \} $ and of the
values
of
$\sigma_n^e \equiv
\sum\limits_{i=1}^\infty \left( \frac 1{z_i}\right) ^n$:}

By taking the ordinary natural logarithm of
\begin{equation}
e_q(z)=\prod_{i=1}^\infty \left( 1-\frac z{z_i}\right) ,
\end{equation}
we obtain
\begin{equation}
\begin{array}{c}
\ln \left\{ e_q(z)\right\} =\sum\limits_{i=1}^\infty \ln
\left\{ 1-
\frac
z{z_i}\right\}  \\
=-z\left\{ \sum\limits_{i=1}^\infty \left( \frac
1{z_i}\right)
\right\}
-
\frac{z^2}2\left\{ \sum\limits_{i=1}^\infty \left( \frac
1{z_i}\right)
^2\right\} -%
\frac{z^3}3\left\{ \sum\limits_{i=1}^\infty \left( \frac
1{z_i}\right)
^3\right\}
\ldots  \\ =b(z)
\end{array}
\end{equation}
where the function
\begin{equation}
\begin{array}{c}
b(z)\equiv \sum\limits_{i=1}^\infty \ln \left\{ 1-\frac
z{z_i}\right\}  \\
=-\sum\limits_{n=1}^\infty \frac 1n\sigma _n^ez^n,|z|<|z_1|
\end{array}
\end{equation}
Fig. 7 of Ref. [5] shows the polar part $\rho _i=|z_i|$ of
the first 8 zeros
of  $e_q(z)$ for $\approx 0.1<q<\approx 0.95$.  Note that
$\rho
_i>\rho
_{i-1}\geq \rho _1$ where $\rho _1$ is the modulus of the
first zero. The function $b(z)=\ln \left\{ e_q(z)\right\} $
is thereby
expressed in terms
of the sum rules for the zeros of $e_q(z)$ since
\begin{equation}
\sigma _n^e\equiv \sum_{i=1}^\infty \left( \frac
1{z_i}\right)
^n;n=1,2,\ldots
\end{equation}
By Eq.(20), the multi-sheet Riemann surface of $b(z)= \ln
\{e_q(z) \} $ consists
of logarithmic branch points at the zeros, $z_i$, of
$e_q(z)$.

The basic properties of $e_q(x)$ displayed in Fig. 1 for $q=
0.1$ follow simply from these expressions for $b(u)$. For
instance, the zeros of $e_q(x)$ correspond to where $b(u)$
diverges.  A sign
change of $e_q(x)$ is due to the principal-value phase
change of ``$+i \pi$" at the branch point of $\ln
\left\{ 1- \frac
z{z_i}\right\}$.

Next, to evaluate these sum rules we proceed as in the above
derivation of the power series
representation for $\ln
_q(1+w)$.  We simply expand both sides of
\begin{equation}
\begin{array}{c}
e_q(z)=e^{b(z)} \\
1+\frac z{[1]!}+\frac{z^2}{[2]!}+\ldots =1+\frac
b{1!}+\frac{b^2}{2!}+\ldots
\end{array}
\end{equation}
Equating coefficients then gives a recursive
formula\footnote{These $\sigma_n^e$ sum rules can also be
evaluated [5] by expanding both sides of an infinite-product
representation of $e_q(z)$.  In this way, from $\sigma_n^e$
for
the first few $n$, we first discovered the general formula
Eq.(23) and Eq.(25).  Eq.(23) describes a pattern similar to
that occurring in the reversion (inversion) of power series.}
for
these
sum
rules:
\begin{equation}
\begin{array}{c}
\sigma _1^e=-1 \\
\sigma _n^e=n\left\{ \sum\limits_{l=2}^n\frac{(-
)^l}{l!}\left(
\sum\limits_{\left(
k_1,k_2,\cdots k_l\right) }\frac{\sigma _{k_1}\sigma
_{k_2}\cdots \sigma
_{k_l}}{k_1k_2\cdots k_l}\right) -\frac 1{[n]!}\right\}
,n\geq 2
\end{array}
\end{equation}
The notation in the second summation is explained
following
Eq.(9) for $\ln _q(1+w).$

The first such sum rules are:
\begin{equation}
\begin{array}{c}
\sigma _1^e=-1 \\ \sigma _2^e=1-\frac 2{[2]!} \\ \sigma
_3^e=-
1+\frac
3{[2]!}-\frac 3{[3]!} \\
\sigma _4^e=1-\frac 4{[2]!}+\frac 4{[3]!}-\frac 4{[4]!}+\frac
2{[2]![2]!}
\end{array}
\end{equation}

The values of $\sigma _n^e$ can also
be
directly
obtained from
\begin{equation}
\sigma _n^e=n\sum_{l=1}^n\frac{(-)^l}l\left\{
\sum_{(k_1,k_{2,}\cdots
k_l)}\frac 1{[k_1]![k_2]!\cdots [k_l]!}\right\} .
\end{equation}
Eq.(25) follows by expanding Eq.(19)
\begin{equation}
\begin{array}{c}
b(z)=-\sum\limits_{n=1}^\infty \frac 1n\sigma _n^ez^n=ln(1+y)
\\
=y-\frac{y^2}2+\frac{y^3}3+...
\end{array}
\end{equation}
where
\begin{equation}
\begin{array}{c}
y=e_q(z)-1 \\
=\frac z{[1]!}+\frac{z^2}{[2]!}+\frac{z^3}{[3]!}+...
\end{array}
\end{equation}
and then equating coefficients of $z^n$.

Equivalently, these formulas can be interpreted as
representations of the reciprocals of the ``bracket''
factorials
in terms of sums of the reciprocals of the zeros of $e_q
(z)$:
\begin{equation}
\begin{array}{c}
\frac 1{[2]!}=\frac 1{2!}-\frac 12\sigma _2^e \\ \frac
1{[3]!}=\frac
1{3!}-\frac 12\sigma _2^e-\frac 13\sigma _3^e \\
\frac 1{[4]!}=\frac 1{4!}-\frac 14\sigma _2^e-\frac 13\sigma
_3^e-\frac
14\sigma _4^e+\frac 18(\sigma _2^e)^2
\end{array}
\end{equation}

The results in this subsection also give $\ln \left\{
E_q(z)\right\} $ for
the analogous $E_q(z)$ for $q>1$ by the substitution
$[n]\rightarrow
[n]_J$.

\subsection{Logarithms and sum rules for related q-analogue
functions:}

(i) For the ``r-$th$'' derivative of
$e_q(x)$, $e_q^{(r)}(x)\equiv \frac{d^r}{dx^r}e_q(x)$, we
similary obtain
[$\alpha
_r\equiv \frac{r!}{[r]!}$]
\begin{equation}
\begin{array}{c}
\ln \left\{ e_q^{(r)}(x)\right\} =\ln \alpha
_r+b^{(r)}(x);r=1,2,\ldots  \\
b^{(r)}(z)=\sum\limits_{i=1}^\infty \ln \left( 1-
\frac
z{z_i^{(r)}}\right)
\end{array}
\end{equation}
where the sum rules for the zeros of the ``r-$th$''
derivative
of $e_q(x)$
are
\begin{equation}
\sigma _n^{(r)}\equiv \sum\limits_{i=1}^\infty \left( \frac
1{z_i^{(r)}}\right) ^n.
\end{equation}

The values of these $e_q (z)$ derivative sum rules are
\begin{equation}
\begin{array}{c}
\sigma _1^{(r)}=-
\frac{r+1}{[r+1]} \\ \sigma _n^{(r)}=n\left\{
\sum\limits_{l=2}^n\frac{(-)^l}{l!}%
\left( \sum\limits_{(k_1,k_2,\cdots k_l)}\frac{\sigma
_{k_1}^{(r)}\sigma
_{k_2}^{(r)}\cdots \sigma _{k_l}^{(r)}}{k_1k_2\cdots
k_l}\right)
-L_n^{(r)}\right\}
\end{array}
\end{equation}
where the $L_n^{(r)}$term is given by

 \begin{equation}
\begin{array}{c}
L_n^{(r)}=
\frac{(n+r)(n+r-1)\cdots (n+1)}{[n+r]!}\frac 1{\alpha _r} \\
=\frac{%
(r+n)(r+n-1)\cdots (r+1)}{[r+n][r+n-1]\cdots [r+1]}\frac
1{n!}
\end{array}
\end{equation}
Equivalently,
\begin{equation}
\sigma _n^{(r)}=n\sum_{l=1}^n\frac{(-)^l}l\left\{
\sum_{(k_1,k_2,\cdots
k_l)}L_{k_1}^{(r)}L_{k_2}^{(r)}\cdots L_{k_l}^{(r)}\right\}
\end{equation}

Thus, the ``r-$th$'' derivative  of $e_q(z)$ is
\begin{equation}
e_q^{(r)}(z)=\frac{r!}{[r]!}\exp \left\{ b^{(r)}(z)\right\}
\end{equation}
where $b^{(r)}(z)=-\sum\limits_{n=1}^\infty \frac 1n\sigma
_n^{(r)}z^n$, $|z|<|z_1^{(r)}|$.

(ii) For the ``r-$th$'' integral of $e_q(z)$ which is defined
in Eq.(15), we obtain
[$\beta
_r\equiv \frac
1{r!}$]

\begin{equation}
\begin{array}{c}
\ln \left\{
\frac{e_q^{(-r)}(x)}{x^r}\right\} =\ln \beta _r+b^{(-
r)}(x);r=1,2,\ldots  \\
b^{(-r)}(z)=\sum\limits_{i=1}^\infty \ln \left( 1-\frac
z{z_i^{(-
r)}}\right)
\end{array}
\end{equation}
where the associated sum rules are
\begin{equation}
\sigma _n^{(-r)}\equiv \sum\limits_{i=1}^\infty \left( \frac
1{z_i^{(-r)}}\right)
^n.
\end{equation}
The values of these $e_q (z)$ integral sum rules are
\begin{equation}
\begin{array}{c}
\sigma _1^{(-r)}=-\frac 1{r+1} \\
\sigma _n^{(-r)}=n\left\{ \sum\limits_{l=2}^n\frac{(-
)^l}{l!}\left(
\sum\limits_{(k_1,k_2,\cdots k_l)}\frac{\sigma _{k_1}^{(-
r)}\sigma
_{k_2}^{(-r)}\cdots \sigma _{k_l}^{(-r)}}{k_1k_2\cdots
k_l}\right) -\frac{r! n!%
}{(r+n)![n]!}\right\}
\end{array}
\end{equation}
Equivalently,
\begin{equation}
\sigma _n^{(-r)}=n\sum_{l=1}^n\frac{(-)^l}l\left\{
\sum_{(k_1,k_2,\cdots
k_l)}L_{k_1}^{(-r)}L_{k_2}^{(-r)}\cdots L_{k_l}^{(-
r)}\right\}
\end{equation}
where the $L_m^{(-r)}$ expression
\begin{equation}
L_m^{(-r)}\equiv \frac{r!m!}{(r+m)![m]!}
\end{equation}
is also the $l=1$ term in Eq.(37).

Thus, the ``r-$th$'' integral of $e_q(z)$ is
\begin{equation}
e_q^{(-r)}(z)=\frac{z^r}{r!}\exp \left\{ b^{(-r)}(z)\right\}
\end{equation}
where $b^{(-r)}(z)=-\sum_{n=1}^\infty \frac 1n\sigma _n^{(-
r)}z^n$, $|z|<|z_1^{(-r)}|$.

(iii) For the q-trigonometric functions, we obtain for the
$cos_q (z)$ function the representation
\begin{equation}
\begin{array}{c}
cos_q(z)=\exp \left\{ b^c(z)\right\}  \\
b^c(z)=\sum\limits_{i=1}^\infty \ln \left( 1-\left(\frac
z{c_i}\right) ^2\right) =-\sum\limits_{n=1}^\infty \frac
1n\sigma _{2n}^cz^{2n}, |z|<|c_1|
\end{array}
\end{equation}
where
\begin{equation}
\sigma _{2n}^c\equiv \sum\limits_{i=1}^\infty \left( \frac
1{c_i^2}\right) ^n.
\end{equation}
The values of the cosine sum rules are
\begin{equation}
\begin{array}{c}
\sigma _2^c=\sum\limits_{i=1}^\infty \left( \frac
1{c_i}\right)
^2=\frac 1{[2]!} \\
\sigma _4^c=\sum\limits_{i=1}^\infty \left( \frac
1{c_i}\right)
^4=\left( \frac
1{[2]!}\right) ^2-\frac 2{[4]!} \\
\sigma _6^c=\sum\limits_{i=1}^\infty \left( \frac
1{c_i}\right)
^6=\left( \frac
1{[2]!}\right) ^3-\frac 3{[2]![4]!}+\frac 3{[6]!} \\
\sigma _{2n}^c=n\left\{ \sum\limits_{l=2}^n\frac{(-
)^l}{l!}\left(
\sum\limits_{(k_1,k_2,\cdots k_l)}\frac{\sigma
_{2k_1}^c\sigma
_{2k_2}^c\cdots
\sigma _{2k_l}^c}{k_1k_2\cdots k_l}\right) -\frac{(-
)^n}{[2n]!}\right\}
\end{array}
\end{equation}
Equivalently,

\begin{equation}
\sigma _{2n}^c=n\sum_{l=1}^n\frac{(-)^l}l\left\{
\sum_{(k_1,k_2,\cdots
k_l)}L_{2k_1}^cL_{2k_2}^c\cdots L_{2k_l}^c\right\}
\end{equation}
where as in Eq.(43)
\begin{equation}
L_{2m}^c\equiv \frac{(-)^m}{[2m]!}
\end{equation}

For the $\sin _q(z)$ function, we find
\begin{equation}
\begin{array}{c}
\sin _q(z)=z\exp \left\{ b^s(z)\right\}  \\
b^s(z)=\sum\limits_{i=1}^\infty \ln \left( 1-\left(\frac
z{s_i}\right) ^2\right)=-\sum\limits_{n=1}^\infty \frac
1n\sigma
_{2n+1}^sz^{2n}, |z|<|s_1|
\end{array}
\end{equation}
where
\begin{equation}
\sigma _{2n+1}^s\equiv \sum\limits_{i=1}^\infty \left( \frac
1{s_i^2}\right) ^n.
\end{equation}
The values of these sine sum rules are
\begin{equation}
\begin{array}{c}
\sigma _3^s=\sum\limits_{i=1}^\infty \left( \frac
1{s_i}\right)
^2=\frac 1{[3]!} \\
\sigma _5^s=\sum\limits_{i=1}^\infty \left( \frac
1{s_i}\right)
^4=\left( \frac
1{[3]!}\right) ^2-\frac 2{[5]!} \\
\sigma _7^s=\sum\limits_{i=1}^\infty \left( \frac
1{s_i}\right)
^6=\left( \frac
1{[3]!}\right) ^3-\frac 3{[3]![5]!}+\frac 3{[7]!} \\
\sigma _{2n+1}^s=n\left\{ \sum\limits_{l=2}^n\frac{(-
)^l}{l!}\left(
\sum\limits_{(k_1,k_2,\cdots k_l)}\frac{\sigma
_{2k_1+1}^s\sigma
_{2k_2+1}^s\cdots
\sigma _{2k_l+1}^s}{k_1k_2\cdots k_l}\right) -\frac{(-
)^n}{[2n+1]!}\right\}
\end{array}
\end{equation}
Equivalently,
\begin{equation}
\sigma _{2n+1}^s=n\sum_{l=1}^n\frac{(-)^l}l\left\{
\sum_{(k_1,k_2,\cdots
k_l)}L_{2k_1+1}^sL_{2k_2+1}^s\cdots L_{2k_l+1}^s\right\}
\end{equation}
where as in Eq.(48)
\begin{equation}
L_{2m+1}^s\equiv \frac{(-)^m}{[2m+1]!}
\end{equation}

\section{Concluding Remarks:}

(1) The above sum rules and logarithmic results are
representation
independent; i.e. they also hold for Jackson's q-exponential
function $E_q(z)
$, its derivatives, integrals, and as well for  its
associated trigonometic
functions $Cos_q(z)$ and $Sin_q(z)$. The only
change is that the
bracket, or deformed integer, $[n]$ is to be replaced by
[$n]_J\equiv \frac{%
1-q^n}{1-q}$.

Since [7,5] the zeros of $E_q(z)$ for $q>1$ are at
\begin{equation}
z_i^E=\frac{q^i}{1-q},
\end{equation}
simple expressions follow: The values of the associated sum rules are
\begin{equation}
\begin{array}{c}
\sigma _n^E\equiv \sum_{i=1}^\infty \left( \frac
1{z_i^{E}}\right) ^n \\
=-
\frac{(1-q)^n}{1-q^n} \\ =-\frac{(1-q)^{n-1}}{[n]_J}.
\end{array}
\end{equation}
A power series representation for the associated natural logarithm is
\begin{equation}
\begin{array}{c}
b^E(z)\equiv \ln \{E_q(z)\} \\
=\sum_{i=1}^\infty
\frac{(1-q)^n}{n(1-q^n)}z^n \\ =\sum_{i=1}^\infty \frac{(1-q)^{n-
1}}{n[n]_J}%
z^n,|z|< | \frac q{1-q} | .
\end{array}
\end{equation}

For both representations, $[n]$ and $[n]_J$, of the derivatives and
integrals of $e_q(z),$ and of the $\cos _q(z)$ and $\sin _q(z)$
functions,
asymptotic formula for their associated zeros are given in Ref.[5] so
simple
expressions also follow for their $\sigma _n$'s and $b(z)$'s in the
regions
where these asymptotic formula apply.

(2) Useful checks on the above results and for use in
applications of
them include: \newline (i) in the bosonic CS(coherent state)
limit
$q\rightarrow 1$,
the normal numerical values must be obtained, \newline (ii)
in the
$q\rightarrow 0$
limit, results corresponding [9] to fermionic CS's should be
obtained [this is a
quick, though quite trivial, check], \newline (iii) by the
use
of
$[n]\rightarrow
$[$n]_J\equiv \frac{1-q^n}{1-q}$, the known exact zeros of
$E_q(z)$ for $q>1$
can be used for non-trivial checks. These zeros are at
$z_i^E=q^n/(1-q)$.

(3) The determination of the series expansion and a general
representation
for the usual natural logarithm for the q-exponential
function, $b(z)=\ln
\{e_q(z) \}$, means that the q-analogue coherent states can
now
be written in
the form of an exponential operator acting on the vacuum
state:
\begin{equation}
\begin{array}{c}
|z\rangle _q=N(|z|)\sum\limits_{n=0}^\infty
\frac{z^n}{\sqrt{[n]!}}|n\rangle _q \\ =N(|z|)\exp \left\{
b(za^{+})\right\}
|0\rangle _q
\end{array}
\end{equation}
where

\begin{equation}
\begin{array}{c}
b(za^{+})=\sum\limits_{i=1}^\infty \ln \left\{ 1-
\frac{za^{+}}{z_i}\right\}  \\ b(za^{+})=za^{+}-\frac
1{[2]!}(za^{+})^2-\left\{ \frac 1{[3]!}-2\left( \frac
1{[2]!}\right)
^2\right\} (za^{+})^3 \\
-\left\{ \frac 1{[4]!}-\frac 5{[3]![2]!}+5\left( \frac
1{[2]!}\right)
^3\right\} (za^{+})^4+\ldots
\end{array}
\end{equation}

(4)The successful evaluations and applications of the sum
rules for the
q-trigonometric functions motivate the following definitions
of q-analogue
generalizations of the usual Bernoulli numbers:
\begin{equation}
\begin{array}{c}
\frac{2^{2n-1}}{(2n)!}B_n^q\equiv \sum\limits_{i=1}^\infty
\left(
\frac
1{s_i}\right) ^{2n} \\ =\sigma _{2n+1}^s
\end{array}
\end{equation}

\begin{equation}
\begin{array}{c}
\frac{2^{2n-1}}{(2n)!}\widetilde{B}_n^q\equiv \frac 1{\left(
2^{2n}-1\right)
}\sum\limits_{i=1}^\infty \left( \frac 1{c_i}\right)
^{2n} \\
=\frac 1{\left( 2^{2n}-1\right) }\sigma _{2n}^c
\end{array}
\end{equation}
Hence, under q-deformation, the usual Bernoulli numbers
become the values of the sum rules for the reciprocals of the
zeros of the q-analogue trigonometric functions, $cos_q (z)$
and $sin_q (z)$.
For the Riemann zeta function, these results do not yield a
unique
definition. However, analogous simple definitions for $p$
complex are
\begin{equation}
\frac 1{\pi ^p}\zeta _q(p)\equiv \sum_{i=1}^\infty \left(
\frac
1{s_i}\right) ^p
\end{equation}
\begin{equation}
\frac 1{\pi ^p}\widetilde{\zeta _q}(p)\equiv \frac 1{\left(
2^p-1\right)
}\sum_{i=1}^\infty \left( \frac 1{c_i}\right) ^p
\end{equation}

``Note added in proof:'' The ordinary natural logarithm of
$E_q(z)$ for
$0<q<1$ is
shown to be related to a q-analogue dilogarithm, $Li_2(z;q)$,
in [10] and
in the
recent survey of q-special functions by Koornwinder [11]:
From
Eq.(53) and $%
E_s(x)E_{1/s}(-x)$ =1, for $0<q<1$%
\begin{equation}
\ln \{E_q(\frac z{1-q})\}=\sum_{i=1}^\infty \frac 1{n(1-
q^n)}z^n \equiv Li_2(z;q)
\end{equation}
which is identical with Eq.(53). Formally [10],
\begin{equation}
\lim _{q\uparrow 1}(1-q)Li_2(z;q)=\sum_{n=1}^\infty
\frac{z^n}{n^2}=Li_2(z)
\end{equation}
the ordinary Euler dilogarithm. Other recent works on q-
exponential
functions are in Refs.[12].

\begin{center}
{\bf Acknowledgments}
\end{center}
This work was
partially
supported by U.S. Dept. of Energy
Contract No. DE-FG 02-96ER40291.

\newpage

\begin{center}
{\bf Figure Captions}
\end{center}

Figure 1: Plot showing the behaviour of the q-analogue
exponential function $ e_q(x) $ for $x$ negative. The $q=0.1$
curve displays the universal behaviour of $ e_q(x) $ for
$q<q_1^*(q_1^*\approx 0.14)$.  As $q$ increases above the
first
collision point at $q_1^*\approx 0.14$, the zeros,
$\mu_i=z_i$, collide in pairs and then move off into the
complex $z$
plane.
They move off as (and remain) a complex conjugate pair.  The
$q=0.2$ curve displays the behaviour of $ e_q(x) $ after the
collision of the first pair of zeros $\mu_1, \mu_2$ but
before the collsion of the first pair of
turning
points.  The first two turning points $\tau_1, \tau_2$
collide
at $ q^*\approx 0.25$.  The turning points $\tau_i$ of
$e_q(z)$ are mapped into the
branch points $b_i$, of $ln_q(w)$.

Figure 2: These two figures and Figs. 3a and 3b show the
Riemann
sheet structure and the mappings of Jackson's exponential
function $E_q(z)$ and of its inverse function $Ln_q(w)$ for
$q^E= 1.09$.  For instance, $w=E_q(z)$ maps the region
labeled $`` 1, 2, 1_L, 2_L "$ \newline in Fig. 2b onto the
upper-half-plane (uhp) of the first $w$ sheet for $Ln_q(w)$,
see Fig. 3a. The turning points $\tau_1,
\tau_2$ are mapped respectively
into the branch points $b_1, b_2$ of Fig. 3a.  These figures
suffice to illustrate the
behaviour of $E_q(z)$ and $Ln_q(w)$ for all $q^E>1$ because
as
$q^E
\rightarrow
1$, the zeros and turning points of
$E_q(z)$ do not collide, but simply move along the negative
$x$ axis and out to
infinity.  In the complex $w$
plane
the associated branch points of $Ln_q(w)$ all move into the
origin.  This
limit thereby gives the usual Riemann surface for $exp(z)$
and $ln(w)$.  Figs. 2 and 3 also illustrate the
Riemann surface for $ e_q(z) $ and
$ln_q(w)$ but only prior to
the collision of the first pair of zeros, i.e. for
$q<q_1^*(q_1^*\approx 0.14)$. Figures 4-8 show the Riemann
surfaces of $ e_q(z) $ and
$ln_q(w)$ for larger $q$ values,$q_1^* < q \leq 1$.

Figure 3: (a) The first upper sheet of $Ln_q(w)$ for $q^E =
1.09$.  The turning points
$\tau_1, \tau_2$ in Fig. 2 for
$E_q(z)$ are mapped respectively into the square-root branch
points $b_1, b_2$ of Fig. 3a, 3b for $Ln_q(w)$. An ``opening
spiral", instead of the usual unit circle, is the ``image"
of the positive $y$ axis (the $x=0$ line) in Fig. 2. The
first lower sheet of $Ln_q(w)$ is the mirror
image of this figure (the reflection is thru the horizontal
$u$ axis); the lower sheets corresponding to the other
``upper sheet" figures in this paper are similarly obtained.
(b) The second upper sheet of $Ln_q(w)$ for $q^E = 1.09$.
Note that the opening spiral continues that in (a).
The cut above the real axis from $b_2$ to $\infty$ goes back
down to the first sheet, Fig. 3a.

Figure 4:  This figure and Fig. 5 show the Riemann
sheet structure and the mappings of $e_q(z)$ and of its
inverse function $ln_q(w)$ for
$0.14 < q \approx 0.22 < 0.25$.  For this range of $q$, the
first two zeros $\mu_1, \mu_2$ of $e_q(x)$ have collided and
have moved off as a complex conjugate pair $\mu_A,
\mu_{\bar{A}}$; the $\mu_A$ zero is marked in this figure.
Note that as in Fig. 2, $Im \lbrace e_q(z) \rbrace =0$ on all
``solid" contour lines, whereas $Re \lbrace e_q(z) \rbrace
=0$
on all ``dashed" coutour lines.

Figure 5:  The first upper sheet for $ln_q(w)$ for $0.14 <
q\approx 0.22 < 0.25$.  When q is increased to
$q\approx0.25$,
the branch points $b_1=b_2$ coincide since the turning points
$\tau_1,
\tau_2$ of Fig. 4 have collided.  Then, the branch cut to the
first
lower sheet nolonger exists.  $\tau_1, \tau_2$ become
a complex conjugate pair $\tau_A, \tau_{\bar{A}}$ and move
off into the complex $z$ plane, as shown in Figs. 6-8.

Figure 6: This figure and Figs. 7-8 show the Riemann
sheet structure and the mappings of $e_q(z)$ and of its
inverse function $ln_q(w)$ for
$q\approx 0.35$. The first two turning points $\tau_1,
\tau_2$
of $e_q(x)$ have collided and have moved off as a complex
conjugate pair $\tau_A, \tau_{\bar{A}}$; the $\tau_A$ turning
point is marked in this figure, $\tau_A= -3.54 + i 1.33$.
The line corresponding to the $\alpha^{'} \beta^{'}$ branch
cut
thru $b_A$, see Figs. 7-8, is the wiggly line from $\alpha$
on the $x<0$ axis, thru $\tau_A$, and on to $\beta$ on the
$Im\lbrace
e_q(z) \rbrace = 0$ curve.  $\tau_A$(and $b_A$) are fixed,
but $\alpha$ and $\beta$( $\alpha^{'}$ and $\beta^{'}$) are
simple
though arbitrary positions on their respective $Im \lbrace
e_q(z) \rbrace =0$ lines.  The third zero $\mu_3$ of $e_q(z)$
is still on the $x<0$ axis.

Figure 7: (a) The first upper sheet of $ln_q(w)$ for $q =
0.35$.  The image of the $x=0$ line in the complex $z$ plane
is shown.  (b) An enlargement of the first quadrant which
shows the $\alpha^{'} \beta^{'}$ branch cut.  For clarity of
illustration, the position of $b_A$ has been displaced from
its true position at $b_A = 0.0222 + i 0.0188$.

Figure 8:  The second upper sheet of $ln_q(w)$ for $q =
0.35$.  The $b_A$ square-root branch point only occurs on the
first two upper sheets, i.e. in Fig. 7  and here.  The
$\alpha^{'}$ point (not shown) lies opposite the $\beta^{'}$
point
and
to the left of the $b_A$ cut structure.

{\bf FOR HARD COPY OF FIGURES, SEND EMAIL TO \newline
cnelson@bingvmb.cc.binghamton.edu }
\end{document}